\documentclass[leqno,12pt]{article}
\usepackage{latexsym}
\usepackage{graphicx}
\usepackage{amsmath}
\usepackage{amssymb}
\usepackage{mathtools}
\usepackage{cite}
\usepackage{color}
\usepackage{subfig}
\usepackage[margin=1.3in, top=1.3in, bottom=1.3in]{geometry}
\usepackage{algorithm}
\usepackage{algorithmic}
\usepackage{hyperref}

\newcommand{\nc}{\newcommand}
\nc{\nt}{\newtheorem}
\nt{thm}{Theorem}[section]
\nt{cor}[thm]{Corollary}
\nt{prop}[thm]{Proposition}
\nt{lem}[thm]{Lemma}
\nt{defn}[thm]{Definition}
\nt{rem}[thm]{Remark}
\nt{exa}[thm]{Example}
\nt{ass}[thm]{Assumption}
\nt{alg}[thm]{Algorithm}
\nt{conj}[thm]{Conjecture}
\nt{claim}[thm]{Claim}
\nt{oracle}[thm]{Oracle}
\nc{\ip}[2]{\mbox{$\langle #1,#2 \rangle$}}
\nc{\pf}{\noindent{\bf Proof\ \ }}
\nc{\finpf}{\hfill{$\Box$}\linespace}
\nc{\linespace}{\vspace
{\baselineskip} \noindent}
\nc{\R}{{\mathbf R}}
\nc{\X}{{\mathbf X}}
\nc{\Y}{{\mathbf Y}}
\nc{\E}{{\mathbf E}}
\nc{\B}{{\mathbf B}}
\nc{\Sn}{{\mathbf S}}
\nc{\Hn}{{\mathbf H}}
\nc{\oR}{\overline{\R}}
\nc{\M}{\mathcal M}
\nc{\e}{\epsilon}

\nc{\Rn}{{\mathbf R}^n}
\nc{\inT}{\mbox{\rm int}\,}
\nc{\cl}{\mbox{\rm cl}\,}

\def\tto{\;{\lower 1pt \hbox{$\rightarrow$}}\kern -12pt
           \hbox{\raise 2.8pt \hbox{$\rightarrow$}}\;}
\newenvironment{myequation}{\setcounter{equation}{\value{thm}}
   \begin{equation}}{\addtocounter{thm}{1}\end{equation}}

\nc{\bmye}{\begin{myequation}}
\nc{\emye}{\end{myequation}}

\begin{document}
\title{
Horoballs and the subgradient method
}
\author{
A.S. Lewis
\thanks{ORIE, Cornell University, Ithaca, NY.
\texttt{people.orie.cornell.edu/aslewis} 
\hspace{2cm} \mbox{~}
Research supported in part by National Science Foundation Grant DMS-2006990.}
\and
G. L\'opez-Acedo
\thanks{Department of Mathematical Analysis -- IMUS, University of Seville, 41012 Seville, Spain
\texttt{glopez@us.es} 
}
\and
A. Nicolae
\thanks{Department of Mathematics, Babe\c{s}-Bolyai University, 400084, Cluj-Napoca, Romania
\texttt{anicolae@math.ubbcluj.ro} 
}
}
\date{\today}
\maketitle

\begin{abstract}
To explore convex optimization on Hadamard spaces, we consider an iteration in the style of a subgradient algorithm.  Traditionally, such methods assume that the underlying spaces are manifolds and that the objectives are geodesically convex:  the methods are described using tangent spaces and exponential maps.  By contrast, our iteration applies in a general Hadamard space, is framed in the underlying space itself, and relies instead on horospherical convexity of the objective level sets.  For this restricted class of objectives, we prove a complexity result of the usual form.  Notably, the complexity does not depend on a lower bound on the space curvature.  We illustrate our subgradient algorithm on the minimal enclosing ball problem in Hadamard spaces.
\end{abstract}
\medskip

\noindent{\bf Key words:} convex, subgradient method, horoball, Hadamard space, complexity, minimal enclosing ball
\medskip

\noindent{\bf AMS Subject Classification:}  90C48, 65Y20, 49M29

\section{Introduction}
We consider optimization problems posed in a geodesic metric space $\M$.  Specifically, we consider first-order algorithms for minimizing an objective function \mbox{$f \colon \M \to \R$} over a set 
$X \subset \M$.  In non-Euclidean spaces $\M$, such algorithms are more subtle than their Euclidean counterparts.  In particular, the literature on subgradient methods in Riemannian manifolds $\M$ generally assumes that $f$ is geodesically convex, and in contrast with the classical Euclidean case, frames algorithms using local linearization, relying on the tangent space to $\M$ at each iterate.  For a thorough exposition, see \cite{boumal2022intromanifolds}.  

This standard approach to subgradient methods on manifolds, originating with \cite{ferreira-oliveira}, thus has two related disadvantages:  the technicalities inherent in local linearization, and the resulting intuitive and formal challenge in extensions to more general spaces.  Optimization in general Hadamard spaces has fascinating applications, including phylogenetic trees \cite{billera} and robotics \cite{ardila-society}.  Subgradient methods are an inviting possibility in such settings:  existing algorithms all rely on proximal techniques \cite{bacak-convex}, which despite their elegant theory, are only rarely implementable.

Motivated both by conceptual simplicity and by generality, therefore, our current development explores an alternative approach to understanding subgradient methods.  The iteration we propose makes no reference to tangent spaces or exponential maps, relying instead on a notion of convexity in Hadamard spaces distinct from geodesic convexity.  This idea --- horospherical convexity --- has its roots in hyperbolic geometry \cite{santalo-yanez,borisenko-miquel}.  

We first diverge from the standard approach to Riemannian subgradient methods in considering objectives that are quasiconvex rather than convex:  in other words, we consider objectives that may not be convex, but have convex lower level sets.  In the Euclidean case, this generalization originated with \cite{kiwiel-quasi}.  The generalization aside, our principal aim here is to focus on the geometry of level sets.  Our objectives are continuous, and typically coercive, so the level sets we consider are usually compact, with nonempty interior.

In Euclidean space, the fact that geodesic and horospherical convexity coincide is the central tool of convex analysis:  the supporting hyperplane theorem.  By contrast, in the case of hyperbolic space, for compact sets with nonempty interior, horospherical convexity entails a certain level of curvature in the boundary, and is a strictly stronger property than geodesic convexity.  The forward implication follows from \cite[Theorem 4.2]{hirai}, for example.  Conversely, in the Poincar\'e disk, where geodesics are segments of circles intersecting the unit circle orthogonally and horospheres are circles tangent to the unit circle, the upper half-disk of radius $\frac{1}{2}$ is geodesically convex but not horospherically convex, not being contained in any horosphere passing through zero.  Nonetheless, applications in Hadamard spaces often involve sets that are both geodesically and horospherically convex, such as intersections of metric balls.   Several recent works have noted the potential applications of horospherical ideas for non-Euclidean convex optimization, both in theory \cite{hirai,bento-fenchel} and practice \cite{fan-yang-vemuri}.  However, the subgradient method that we present here appears new.

Complexity analysis for standard projected subgradient methods on Hadamard manifolds originated with the seminal work \cite{zhang-sra}, which shows in particular that the mean excess objective value over $n$ iterations behaves like $O(\frac{1}{\sqrt{n}})$,  just as in the Euclidean version.  Interestingly, the hidden constant deteriorates as the curvature becomes more negative, as subsequent analysis \cite[Theorem 8]{criscitiello-boumal} shows that it must, at least for general geodesically convex objectives.

Our second and principal point of divergence from \cite{zhang-sra} and other standard non-Euclidean subgradient developments lies in our avoidance of any reference to tangent spaces and exponential maps.  Rather than situating the algorithm in a manifold, our subgradient-style method makes sense in any Hadamard space.  Our complexity analysis is geometrically appealing:  like the algorithm itself, it makes no use of tangent spaces, and as we have noted, applies to objectives that may be nonconvex but have horospherically convex level sets.  Strikingly, unlike \cite{zhang-sra}, our complexity result and proof involves no lower curvature bound: the lower bound on complexity in \cite[Theorem 8]{criscitiello-boumal} does not apply, because we only consider objectives with horospherically convex level sets. 

To summarize our contributions, we present a new horospherical approach to subgradient methods that we believe is geometrically illuminating.  Our algorithm is framed in a general Hadamard space, and relies on no tangent constructions.  The algorithm does not apply to every geodesically convex objective, but its range of application is still rich.  Furthermore, our complexity result involves no lower curvature bound, so even applies in spaces with infinite negative curvature and bifurcating geodesics, such as the tree spaces and cubical complexes of \cite{billera,ardila-society}.  We apply our horospherical subgradient algorithm in particular to the minimal enclosing ball problem.  In Euclidean space, and more generally in Hadamard manifolds, our algorithm amounts to the approach of \cite{badoiu-clarkson,arnaudon-nielsen}, but algorithms for the minimal enclosing ball problem in a general Hadamard space were not previously available.  We present a small computational illustration in a CAT(0) complex.

\section{Examples}
In a Hadamard space $(\M,d)$, many interesting optimization problems involve simple combinations of distance functions.  To ground our discussion, we consider two representative examples.

\begin{exa}[Circumcenters] \label{circumcenters}
{\rm
Any nonempty finite set $A \subset \M$ has a unique {\em circumcenter} \cite[Proposition II.2.7]{bridson}, which is the unique minimizer $\bar x$ of the function $f \colon \M \to \R$ defined by
\[
f(x) ~=~ \max_{a \in A} d(x,a) \qquad (x \in \M).
\]
The value $f(\bar x)$ is called the {\em circumradius}.  Other names for this optimization problem include the {\em smallest} or {\em minimal enclosing ball problem}, and the {\em $1$-center problem}.  For a survey, see \cite{arnaudon-nielsen}.
}
\end{exa}

\begin{exa}[Intersecting balls] \label{balls}
{\rm
Consider any nonempty finite set $A \subset \M$.  Given a radius $\rho_a \ge 0$ for each point $a \in A$, the corresponding balls $B_{\rho_a}(a)$ have nonempty intersection if and only if the function $f \colon \M \to \R$ defined by
\[
f(x) ~=~ \max_{a \in A} \{d(x,a) - \rho_a \} \qquad (x \in \M).
\]
has nonpositive minimum value.
}
\end{exa}

Both examples involve minimizing an objective function belonging to the following broad class.

\begin{defn} \label{envelope}
{\rm
A function $f \colon \M \to \R$ is a {\em distance envelope} if there exists a compact topological parameter space $\Psi$, and continuous functions
\[
a \colon \Psi \to \M, \qquad \beta \colon \Psi \to \R_+, \qquad \gamma \colon \Psi \to \R, 
\]
such that
\bmye \label{max}
f(x) ~=~ \max_{\psi \in \Psi} \{ \beta(\psi) d\big(x,a(\psi)\big) + \gamma(\psi) \} \qquad (x \in \M).
\emye
The points $a(\psi)$ are the {\em centers}.
}
\end{defn}

\section{Horoballs}
The subgradient-style method we develop, rather than using the local idea of the tangent space, instead relies on global constructions fundamental in CAT(0) geometry, namely horoballs.  We consider a Hadamard space $(\M,d)$ with the {\em geodesic extension property\/}:  in other words, for any geodesic segment 
$[x,y]$ in $\X$, there exists a {\em ray} $r \colon \R_+ \to \M$ (by which we always mean a unit-speed geodesic ray) such that $r(0) = x$ and $r(t) = y$ for some $t \ge 0$.  Following standard terminology \cite{bridson}, to any ray $r$ we associate the {\em Busemann function} $b_r \colon \M \to \R$ defined by
\[
b_r(x) ~=~ \lim_{t \to \infty}(d\big(x,r(t)\big) - t) \qquad (x \in \M).
\]

Busemann functions are $1$-Lipschitz and convex, and their lower level sets are called {\em horoballs}.  Following standard terminology in hyperbolic geometry \cite[Definition 4.5]{gallego-solanes-teufel}, we call a closed set $F \subset \M$ {\em horospherically convex}  if every boundary point $x$ of $F$ has a {\em supporting ray}, by which we mean a ray $r$ satisfying $r(0) = x$ and  
\[
F ~\subset~ \{z \in \M : b_r(z) \le 0\}. 
\]
In that case we describe the right-hand side as a {\em supporting horoball} for $F$ at $x$.

\begin{exa}[The Euclidean case] \label{euclidean1}
{\rm
Consider the space $\M = \Rn$ with Euclidean distance.  Given any ray $r(t) = x - tu$ (for $t \ge 0$), for a point $x \in \M$ and a unit vector $u \in \Rn$, the corresponding Busemann function is given by
\[
b_r(z) ~=~ \ip{u}{z-x} \qquad (z \in \M).  
\]
Thus horoballs in $\Rn$ are halfspaces.  The ray $r$ supports a closed convex set $F \subset \Rn$ at $x$ if and only if $u$ lies in the normal cone $N_F(x)$, in which case the halfspace $\{z : \ip{u}{z-x} \le 0 \}$ is a supporting horoball for $F$ at $x$.  In particular, 
all closed convex sets $F$ are horospherically convex. 
}
\end{exa}

\begin{defn} \label{polyhedral-def}
{\rm
A function $f \colon \M \to \R$ is {\em horospherically polyhedral} if there exist rays 
$r_i \colon \R_+ \to \M$, for $i=1,2,\ldots,m$, such that
\bmye \label{polyhedral}
f(x) ~=~ \max_{i=1,\ldots,m} \{\beta_i b_{r_i}(x) + \gamma_i\} \qquad (x \in \M).
\emye
}
\end{defn}

\section{The oracles}
The iteration that we propose relies on two computational resources, described next.  The first is standard.

\begin{ass}[Projection oracle] \label{projection-ass}
The feasible region $X \subset \M$ is nonempty, closed, bounded, and geodesically convex.  For any input $x \in \M$, the projection oracle outputs the nearest point in $X$, denoted $\mbox{\rm Proj}(x)$.
\end{ass}

\noindent
The nearest point always exists and is unique \cite[Proposition II.2.4]{bridson}.
The simplest example of the projection oracle is when the feasible region $X$ is just a metric ball $B_\rho(a)$ for some radius $\rho \ge 0$ and some center $a \in \M$.  In that case, the projection 
$\mbox{\rm Proj}(x)$ equals $x$ whenever $d(x,a) \le \rho$, and otherwise is the unique point on the geodesic segment $[a,x]$ at a distance $\rho$ from $a$.

The second assumption is less standard.

\begin{ass}[Support oracle] \label{support-ass}
The objective function $f \colon \M \to \R$ is Lipschitz, with horospherically convex lower level sets
\[
f_x ~=~ \{z  \in \M : f(z) \le f(x)\} \qquad (x \in \M).
\]
For any input $x \in \M$ and any step length $\epsilon > 0$, the support oracle either responds that $x$ minimizes $f$, or it outputs a point $x_\epsilon$ on a supporting ray for $f_x$ at $x$ such that $d(x,x_\epsilon) = \epsilon$.
\end{ass}

We illustrate the support oracle with some examples.  The first shows that, in the classical case of convex optimization in Euclidean space, the support oracle coincides with the standard version.

\begin{exa}[The Euclidean case]~
{\rm
Consider a continuous convex function $f \colon \Rn \to \R$, and any point $x \in \M$ that does not minimize $f$.  Then valid outputs of the support oracle are those points of the form
\[
x_\epsilon = x - \frac{\epsilon}{|g|}g \qquad \mbox{for}~ g \in \partial f(x).
\]
To see this, we recall Example \ref{euclidean1}, and use the standard relationship between the normal cone to the level set and the subdifferential, $N_{f_x}(x) = \R_+ \partial f(x)$.
}
\end{exa}

\begin{exa}[Distance envelopes] \label{ex-envelope}
{\rm
In Definition \ref{envelope}, consider the distance envelope $f$ given by equation (\ref{max}).
Given any point $x \in \M$, fix any parameter value $\psi$ attaining this maximum.  If $\beta(\psi) = 0$, then $x$ minimizes $f$.  On the other hand, if $\beta(\psi) > 0$, then any ray $r \colon \R_+ \to \M$ satisfying $r(0) = x$ and passing through the center $a(\psi)$ supports the level set $f_x$, because any point $z \in f_x$ satisfies
\[
\beta(\psi) d\big(z,a(\psi)\big) + \gamma(\psi) ~\le~ f(z) ~\le~ f(x) ~=~ \beta(\psi) d\big(x,a(\psi)\big) + \gamma(\psi)
\]
and hence
\[
d\big(z,a(\psi)\big) ~\le~ d\big(x,a(\psi)\big),
\]
so
\[
b_r(z) ~=~ \inf_{t \ge 0} \{d\big(z,r(t)\big) - t\} ~\le~ d\big(z,a(\psi)\big) - d\big(x,a(\psi)\big) ~\le~ 0.
\]
The support oracle can therefore output the point $x_\epsilon = r(\epsilon)$.  In summary, for a distance envelope, to implement the support oracle we simply move a distance $\epsilon$ along a ray towards a center attaining the envelope's value.
}
\end{exa}

\begin{exa}[Horospherically polyhedral functions]
{\rm
In Definition \ref{polyhedral-def}, consider the distance envelope $f$ given by equation (\ref{polyhedral}).  Given any point $x \in \M$, fix any index $i$ attaining this maximum.  If 
$\beta_i = 0$, then $x$ minimizes $f$.  On the other hand, suppose $\beta_i > 0$.  There exists a unique ray $r \colon \R_+ \to \M$ asymptotic to $r_i$ and satisfying $r(0) = x$, by \cite[Proposition II.8.2]{bridson}.  We claim that this ray supports the level set $f_x$ at $x$.  Note that the corresponding Busemann function $b_r$ differs from the Busemann function $b_{r_i}$ only by a constant \cite[Exercises II.8.23]{bridson}.  Any point $z \in f_x$ satisfies 
\[
\beta_i b_{r_i}(z) + \gamma_i ~\le~ f(z) ~\le~ f(x) ~=~  \beta_i b_{r_i}(x) + \gamma_i,
\]
so $b_{r_i}(z) \le b_{r_i}(x)$, and hence $b_r(z) \le b_r(x) = 0$.  The support oracle can therefore output the point $x_\epsilon = r(\epsilon)$.
}
\end{exa}

\subsection*{Discussion}
The availability of the support oracle is, of course, a strong assumption.  As we discussed in the introduction, geodesic and horospherical convexity are distinct notions.  Thus the level sets of a geodesically convex objective function $f$ may not be horospherically convex, and even if they are, traditional subgradients may not correspond to the requisite supporting horoballs.  Using standard Riemannian notation, the traditional oracle returns a vector $g$ in the tangent space $T_x(\M$) satisfying the inequality
\[
f(z) ~\ge~ f(x) + \ip{g}{\mbox{Exp}_x^{-1}(z)}_x \qquad (z \in \M),
\]
which ensures that the level set $L_x$ is contained in the set
\[
\{z \in \M :  \ip{g}{\mbox{Exp}_x^{-1}(z)}_x \le 0\}.
\]
This set contains the point $x$, but it may not correspond to any supporting horoball.

On the other hand, the support oracle is appealing geometrically, since it avoids any reference to tangent spaces and exponential maps.  Furthermore, it relies only on convexity properties of the level sets of the objective $f$, rather than on convexity of $f$ itself.  The applicability of the subgradient method to Euclidean quasiconvex minimization was developed by \cite{kiwiel-quasi}.

\section{A subgradient-style method} \label{alg-subgradient}
In a Hadamard space $(\M,d)$, suppose that Assumptions \ref{projection-ass} and \ref{support-ass} hold, and consider the optimization problem
\[
\inf_{x \in X} f(x).
\]

\begin{alg}[Projected subgradient iteration] \label{iteration} \mbox{}\\
{\rm
Choose a fixed step length $\epsilon > 0$.  \\
Repeatedly update the current iterate $x \in X$ as follows.
\begin{itemize}
\item  
Call the support oracle in Assumption \ref{support-ass} at the point $x$.
\item
If the support oracle recognizes that $x$ attains $\min_\M f$, stop.
\item
Else, update $x \leftarrow x_\epsilon$.
\item
Call the projection oracle in Assumption \ref{projection-ass} at the point $x$.
\item
Update $x \leftarrow \mbox{\rm Proj}(x)$.
\end{itemize}
}
\end{alg}

\begin{exa}[Circumcenters and intersecting balls]~~
{\rm
We return to Examples \ref{circumcenters} and \ref{balls}.
Given an arbitrary point $a_0 \in \M$, if the radius $\rho > 0$ is large enough, then both the circumcenter and intersecting problems are equivalent to minimizing their respective objectives $f$ over the ball $X = B_\rho(a_0)$.  Projecting onto the ball $X$ is easy, so we can implement the subgradient iteration simply:  we simply need to compute (and possibly extend) geodesics between iterates and points in $A$ or $a_0$.
}
\end{exa}

\section{Complexity analysis}
The projected subgradient iteration that we have described is applicable in any Hadamard space with the geodesic extension property.  We make no use of tangent spaces, and the curvature of the space may be unbounded below.  In addition to Hadamard manifolds, the result applies in spaces like CAT(0) cubical complexes, which are not manifolds, and have infinite negative curvature.

In addition to Assumptions \ref{projection-ass} and \ref{support-ass}, we collect up the remaining assumptions.

\begin{ass} \label{remaining-ass}
\mbox{}
\begin{itemize}
\item
The Hadamard space $(\M,d)$ has the geodesic extension property.
\item
The diameter of the feasible region $X \subset \M$ is bounded above by $D > 0$.
\item
The objective function $f \colon \M \to \R$ is globally Lipschitz, with constant $L>0$.
\item
The optimal value $f^* = \inf_X f$ is attained at some point $x^* \in X$.
\end{itemize}
\end{ass}

Under all of our assumptions, we prove that the iterates $x^1,x^2,\ldots$ generated by the subgradient method with a suitable step length, ensures a complexity estimate of the usual kind:
\[
\min_{i=1,\ldots,n} f(x_i) ~=~ \min_X f ~+~ O\Big(\frac{1}{\sqrt{n}}\Big).
\]
The complexity estimate depends only on the constants $L$ and $D$, and not on any lower bound on the curvature.

\begin{thm} \label{main}
Consider the optimization problem $\inf_X f$, and suppose that Assumptions \ref{projection-ass}, \ref{support-ass}, and \ref{remaining-ass} hold.  
For any integer $n>0$, consider Algorithm~\ref{iteration} with step length $\epsilon = \frac{D}{\sqrt{n}}$.
After $n$ iterations, the average of the function values at the iterates exceeds $\min_X f$ by a quantity no larger than
\[
\frac{LD}{\sqrt{n}}.
\]
\end{thm}

\pf
Beginning with any current iterate $x \in X$, one iteration of Algorithm~\ref{iteration}, outputs a point $x_\epsilon$ on a supporting ray $r \colon \R_+ \to \M$ for $f_x$ at $x$ such that 
$d(x,x_\e) = \e$.  We denote the corresponding Busemann function $b_r$ simply by $b$.  Since $x^*$ minimizes $f$, we know $f(x^*) \le f(x)$, and hence $b(x^*) \le 0$.
We have $x_\e = r(\e)$, and more generally, for $t \ge 0$, we denote the point $r(t)$ by $x_t$.
 
For any value $t \ge \e$ we have 
\[
x_\e ~=~ \Big(1-\frac{\e}{t}\Big)x + \frac{\e}{t} x_t.
\]
Since $\M$ is a Hadamard space, nonpositive curvature implies
\[
d^2(x^*,x_\e) 
~\le~ 
\Big(1-\frac{\e}{t}\Big)d^2(x^*,x) + \frac{\e}{t}d^2(x^*,x_t)
- \frac{\e}{t}\Big(1-\frac{\e}{t}\Big) d^2(x,x_t)
\]
so
\begin{eqnarray*}
d^2(x^*,x_\e) - \Big(1-\frac{\e}{t}\Big)d^2(x^*,x) - \e^2
&\le&  
\frac{\e}{t} ( d^2(x^*,x_t) - t^2) \\
&=&  
\e( d(x^*,x_t) - t)\Big( \frac{1}{t}d(x^*,x_t) + 1\Big).
\end{eqnarray*}
As $t \to +\infty$, we know $d(x^*,x_t) - t \to b(x^*)$ by definition, and hence $\frac{1}{t} d(x^*,x_t) \to 1$, from which we deduce
\[
d^2(x^*,x_\e) - d^2(x^*,x) - \e^2 ~\le~ 2\e b(x^*).
\]

There exists a unique ray originating at $x^*$ and asymptotic to the ray $r$.  Extend this ray to a point $\tilde x \in \M$ satisfying
$d(x^*,\tilde x) = 1-b(x^*)$, giving a ray $\tilde r \colon \R_+ \to \M$ asymptotic to $r$ satisfying 
$\tilde r(0) = \tilde x$ and $\tilde r\big(1-b(x^*)\big) = x^*$.  The  corresponding Busemann function 
$\tilde b$ satisfies
\[
\tilde b(x^*) ~=~ \tilde b(\tilde r\big(1-b(x^*)\big)) ~=~ b(x^*) - 1,
\]
and it differs from the Busemann function $b$ only by a constant, so 
$\tilde b = b-1$.  For all $t < 1$ we have
\[
b\big(\tilde r(t)\big) ~=~ \tilde b\big(\tilde r(t)\big) + 1 ~=~ 1-t ~>~ 0.
\]
Since the ray $r$ supports the level set $f_x$ at $x$, we know $b(z) \le 0$ for all points $z \in f_x$.
Hence we have $\tilde r(t) \not\in f_x$, or equivalently, $f\big(\tilde r(t)\big) > f(x)$.  By continuity, the point $\hat x = \tilde r(1)$ satisfies $f(\hat x) \ge f(x)$.  Since $d(x^*,\hat x) = -b(x^*)$,
the Lipschitz condition ensures
\[
f(x) - f(x^*) ~\le~ f(\hat x) - f(x^*) ~\le~ -Lb(x^*).
\]
To summarize, we have proved
\[
d^2(x_\epsilon,x^*) - d^2(x,x^*) ~\le~ \epsilon^2 - \frac{2\epsilon}{L} (f(x) - f^*).
\]
Since $x^* \in X$ and the projection $\mbox{Proj}_X$ is nonexpansive, we deduce
\[
\frac{2}{L} (f(x) - f^*) ~\le~ \frac{1}{\epsilon}\big( d^2(x,x^*) - d^2(x_{\mbox{\scriptsize new}},x^*) \big) 
~+~ \epsilon.
\]

Now suppose we initialize the algorithm at a point $x^1 \in X$, and apply the update $n$ times, generating the points $x^1,x^2,x^3,\ldots,x^n,x^{n+1} \in X$.  Summing the corresponding inequalities shows
\[
\frac{2}{L} \Big( \frac{1}{n}\sum_{i=1}^n f(x^i) - f^* \Big) ~\le~ \frac{D^2}{n\epsilon} + \epsilon.
\]
We deduce
\[
\frac{1}{n}\sum_{i=1}^n f(x^i) ~\le~ f^* + \frac{LD}{\sqrt{n}},
\]
as required.
\finpf

\section{Computing circumcenters in Hadamard space}
We conclude by returning to the circumcenter problem, Example~\ref{circumcenters}).  The problem has a long history dating back to Sylvester (1857) \cite{sylvester}, and remains fundamental in computational geometry.  A simple Euclidean algorithm, with complexity $O(\frac{1}{\sqrt{n}})$ appeared in \cite{badoiu-clarkson}, and can be interpreted as the subgradient method \cite{krivosija-munteanu}.  The Euclidean algorithm was extended to Riemannian manifolds in \cite{arnaudon-nielsen}, and we shall see that a constant-step-length version in fact extends further to arbitrary Hadamard spaces, where it can be interpreted as the horospherical subgradient method, Algorithm \ref{iteration}.

This extension is motivated in part by an interesting class of optimization problems in Hadamard space that arises from averaging finite sets in the BHV space of phylogenetic trees introduced in \cite{billera}.  CAT(0) cubical complexes, of which the BHV space is an example, are Hadamard spaces with curvature unbounded below:  in particular, they are not manifolds.  In addition to BHV space, CAT(0) cubical complexes have also proved useful models in robotics, as surveyed recently in \cite{ardila-society}.  Crucially from an algorithmic perspective, geodesics in CAT(0) cubical complexes are efficiently computable \cite{owen-provan,hayashi}.

Consider the problem of averaging a given nonempty finite set $A$ in a Hadamard space $\M$.  Standard notions of average include the {\em Fr\'echet mean} and {\em medians}, which are minimizers of the functions 
\[
\sum_{a \in A} d^2(x,a) \qquad \mbox{and} \qquad \sum_{a \in A} d(x,a) \qquad (x \in \M)
\]
respectively.  Such points can be approximated via a sequence of geodesic computations in a splitting proximal point algorithm.  The convergence of the algorithm is proved in \cite{bacak-medians}, and illustrated on real data in BHV space in \cite{bacak}. While \cite{bacak-medians} states no complexity result, the algorithm and its convergence proof are based on a Euclidean version \cite{bertsekas-incremental}, whose complexity, like that of Algorithm \ref{iteration}, is $O(\frac{1}{\sqrt{n}})$.

The horospherical subgradient method that we have described here opens an alternative possibility for averaging a finite set $A$.  Instead of medians or the Fr\'echet mean, we might instead approximate the circumcenter, by minimizing the objective function 
$f \colon \M \to \R$ defined, as in Example~\ref{circumcenters}:
\bmye \label{objective}
f(x) ~=~ \max_{a \in A} d(x,a) \qquad (x \in \M).
\emye
As we shall see, implementing the horospherical subgradient method to minimize the objective function $f$ is straightforward.

One curiosity suggesting the circumcenter as an alternative averaging technique in CAT(0) cubical complexes is the {\em stickiness} to which the Fr\'echet mean is prone:  for a short exposition, see \cite{miller}. Consider for example the simplest BHV space, the {\em tripod} consisting of three {\em legs\/}:  copies of $\R_+$ glued together at zero.  Three points, one on each leg at distances $a,b,c>0$ from zero, have mean at zero if and only if
\bmye \label{triangle}
a+b>c, \qquad b+c>a, \qquad c+a>b.  
\emye
The point zero is {\em sticky} because the corresponding set of vectors $(a,b,c)$ has nonempty interior.  More precisely, we can identify each vector $(a,b,c) \in \R^3_{++}$ with a probability measure on the tripod supported on the three corresponding points, with equal weights:  this set of measures inherits the topology of $\R^3$, and any measure in the resulting space that satisfies the inequalities (\ref{triangle}) {\em sticks} to zero, in the terminology of \cite{huckemann}.  The median in this example behaves even more simply:  the unique median is zero for all values $a,b,c$.  Notice, in contrast, that zero is a circumcenter if and only if the maximum component of the vector $(a,b,c)$ is not unique.  Since the set of such vectors has measure zero, the point zero is not sticky as a circumcenter. 

To compute the circumcenter using the horospherical subgradient method, we choose an arbitrary point 
$\bar a \in A$, define the radius $\rho = f(\bar a)$, and consider the feasible region 
$X = B_\rho(\bar a)$.    Notice $A \subset X$.  Since all points $x \not \in X$ satisfy 
$f(x) \ge d(x,\bar a) > \rho$, minimizing $f$ over the whole space $\M$ is equivalent to minimizing $f$ over $X$.  The set $X$ satisfies Assumption \ref{projection-ass}, and the projection oracle is simple: 
\[
\mbox{Proj}(x) ~=~ 
\left\{
\begin{array}{ll}
x & (x \in X) \\
\frac{d(x,a)}{\rho}x + (1 - \frac{d(x,a)}{\rho})\bar a) & (x \not\in X).
\end{array}
\right.
\]
As we shall see, in fact the subgradient algorithm only calls the projection oracle for inputs $x \in X$, so the projection plays no role.

The objective function $f$ is $1$-Lipschitz.  By Example \ref{ex-envelope}, its level sets are horospherically convex, and given any input point $x \in \M$ and step length $\e > 0$, the support oracle chooses a point $a \in A$ attaining $\max_A d(x,\cdot) = f(x)$, and then outputs the point $x_\e$ on the ray originating at $x$ and passing through $a$, at a distance $\e$ from $x$.  Thus Assumption \ref{support-ass} holds.  Furthermore, providing $\e \le f(x)$, the point $x_\e$ is a convex combination of $x$ and $a$.

We assume that the Hadamard space $\M$ has the geodesic extension property, as holds in particular for BHV space.  We define the constant $D = 2\rho$.  Assumption \ref{remaining-ass} then holds.

Now consider Algorithm \ref{iteration} with step length $\e = \frac{D}{\sqrt{n}}$ and number of iterations $n \ge 16$.  Notice $\e \le \frac{\rho}{2}$.  If any point 
$x \in \M$ satisfied $f(x) < \e$, then we deduce the contradiction
\[
d(\bar a, a) ~\le~ d(\bar a,x) + d(x,a) ~<~ \e + \e ~\le~ \rho ~=~ f(\bar a) \qquad \mbox{for all}~ a \in A.
\] 
We therefore have $f(x) \ge \e$ for all $x \in \M$.  Consequently, for any iterate $x$ in the feasible region $X$, the output $x_\e$ of the support oracle also lies in $X$.  By induction, the projection oracle thus plays no role, as we claimed.  

We can summarize the algorithm, in terms of the function (\ref{objective}), as follows.  While the approach is similar to the Euclidean algorithm of \cite{badoiu-clarkson} and its Riemannian extension \cite{arnaudon-nielsen}, we do not know of any previous algorithm for circumcenters in arbitrary Hadamard spaces.

\begin{alg}[Subgradient iteration for circumcenter] \label{alg-circumcenter}
\label{bisection}
{\rm
\begin{algorithmic}
\STATE
\STATE 	{\bf input:} nonempty finite $A \subset \M$, number of iterations $n \ge 16$
\STATE	choose $x \in A$						\hfill \% {\em current center}
\STATE	$x_{\mbox{\scriptsize best}} = x$		\hfill \% {\em best center so far}
\STATE	$\rho = 0$ 
\FOR{$a \in A \setminus \{x\}$ }
\STATE	$\rho = \max\{\rho,d(a,x)\}$ 			\hfill \% {\em radius of current ball enclosing $A$}	
\ENDFOR
\STATE	$f_{\mbox{\scriptsize best}} = \rho$	\hfill \% {\em best radius so far}
\STATE	$\e = \frac{2\rho}{\sqrt{n}}$ 			\hfill \% {\em step length}
\FOR{$\mbox{iteration} = 1,2,\ldots,n$}
\STATE	$\gamma=0$
\FOR{$a \in A$}
\STATE	$\delta = d(a,x)$
\IF{$\delta > \gamma$}
\STATE	$\gamma = \delta$ 						\hfill \% {\em radius of current enclosing ball\ldots}
\STATE	$\hat a = a$ 							\hfill \% {\em \ldots and furthest point in $A$}
\ENDIF
\ENDFOR
\IF{$\gamma < f_{\mbox{\scriptsize best}}$}
\STATE	$f_{\mbox{\scriptsize best}} = \gamma$	\hfill \% {\em update best radius \ldots}
\STATE	$x_{\mbox{\scriptsize best}} = x$		\hfill \% {\em \ldots and best center}
\ENDIF
\STATE	$\lambda = \frac{\e}{\gamma}$
\STATE	$x = (1-\lambda)x + \lambda \hat a$		\hfill \% {\em move towards furthest point}
\ENDFOR
\RETURN	$x_{\mbox{\scriptsize best}}$
\end{algorithmic}
}
\end{alg}

\noindent
If the set $A$ has cardinality $k$, then the initial {\tt for} loop requires $k-1$ distance calculations, and then each of the $n$ iterations requires $k$ more distance computations, one of which results in the calculation of a convex combination of two points.  The diameter $\mu$ of the set $A$ satisfies $\mu \ge \rho$, so the feasible region $X$ has diameter no larger than $2\mu$.  Theorem \ref{main} now implies the following complexity result.

\begin{cor}[Complexity of circumcenters] \label{circumcenter-complexity}
Consider any Hadamard space $\M$ with the geodesic extension property, and any subset $A \subset \M$ of cardinality $k>0$ and diameter $\mu$.  Fix any integer $n \ge 16$.  Then Algorithm \ref{alg-circumcenter} terminates after computing $k(n+1) - 1$ distances and $n$ convex combinations in $\M$, and returns a point whose distance from every point in $A$ exceeds the circumradius of $A$ by no more than $\frac{2\mu}{\sqrt{n}}$.
\end{cor}

\subsubsection*{Example:  circumcenter of three points in an orthant space}
To illustrate Algorithm \ref{alg-subgradient} on a polyhedral complex analogous to the tree spaces of \cite{billera}, we consider the geodesic space consisting of five Euclidean quadrants glued together along their edges to form a cycle.  This orthant space $\M$ is a well-known example of a Hadamard space:  it can be modeled isometrically as the subset $\M'$ of the space $\R^3$ consisting of the union of the quadrants
\[\R_+ \times \R_- \times \{0\}, \quad \R_- \times \R_- \times \{0\}, \quad \R_- \times \R_+ \times \{0\}, \quad \{0\} \times \R_+ \times \R_+, \quad \R_+ \times \{0\} \times \R_+,
\]
with the intrinsic metric induced by Euclidean distance \cite[Figure 9]{ardila-society}.

Computing geodesics in this orthant space $\M$ is easy.  Given any two points $x,y \in \M$, after a cyclic permutation of the quadrants, we can suppose that the corresponding points $x'$ and $y'$ in the model both lie in the plane $\R \times \R \times \{0\}$.  In the model, the geodesic $[x',y']$ is just the Euclidean line segment $[x',y']$ in the plane $\R \times \R \times \{0\}$, unless that line segment intersects the open quadrant $\R_{++} \times \R_{++} \times \{0\}$, in which case the geodesic consists of the union of the two geodesics $[x',0]$ and $[0,y']$.

In the model $\M'$, consider the set $A$ consisting of the three points
\[
(0,0,\sqrt{5}) ~~,~~ (1,-2,0) ~~,~~ (-2,1,0).
\]
The circumcenter of $A$ is the point $(0,0,0)$.  To see this, consider the objective function 
$f(x) = \max_A d(x,\cdot)$ (for $x \in \M'$) whose unique minimizer is the circumcenter $(p,q,r)$ of $A$.  By symmetry, the point $(q,p,r)$ also lies in $\M'$, and the value of $f$ there is identical, so in fact $(q,p,r) = (p,q,r)$, and hence $p=q$.  Consequently the circumcenter of $A$ must lie on the set
\[
\{(0,0,t) : t \ge 0\} ~\cup~ \{(t,t,0) : t \le 0 \}.
\]
A quick calculation shows that $f$ is minimized on this set at the point $(0,0,0)$.

\begin{figure}
    \centering
    \includegraphics[width=0.8\textwidth]{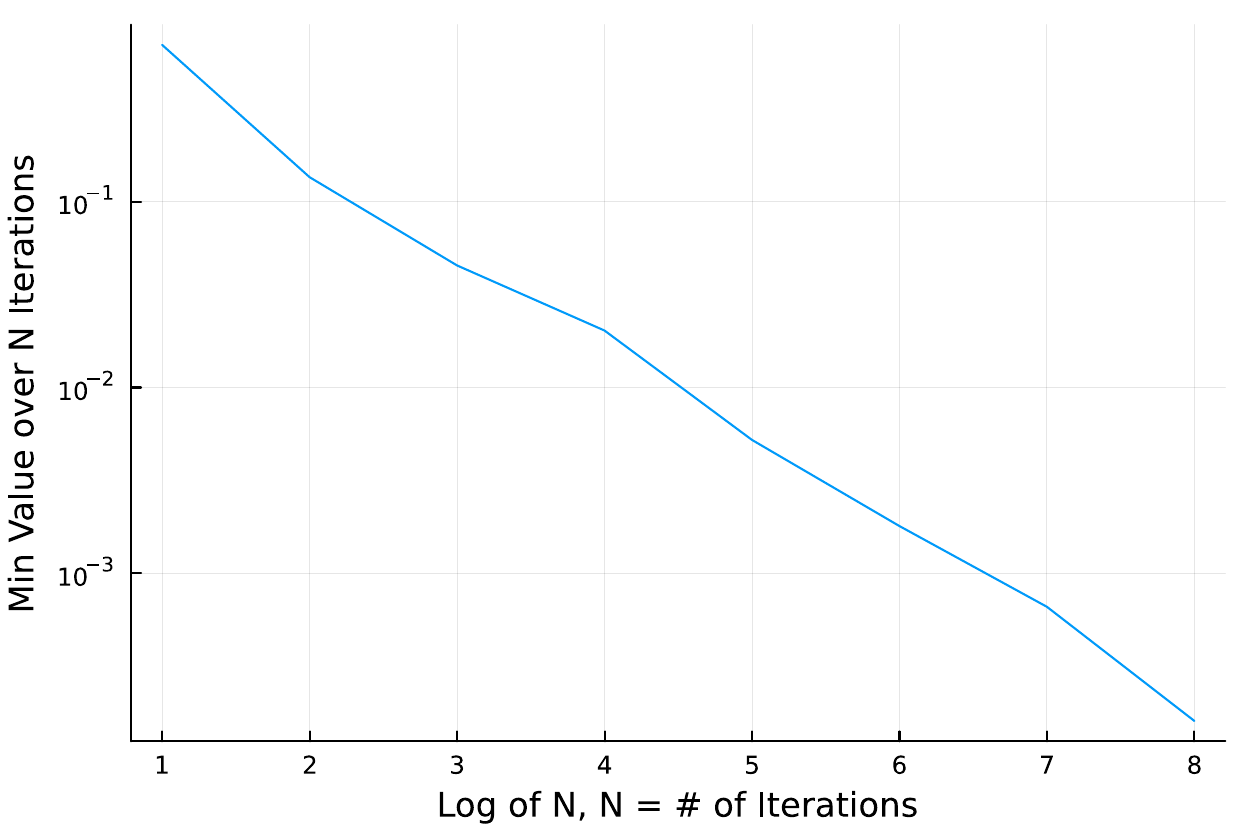}
    \caption{Algorithm \ref{alg-circumcenter} applied to the circumcenter and circumradius $\sigma$ of three points in a Hadamard space of five orthants.  For number of iterations $N=10^1, 10^2,\ldots,10^8$, the plot shows the minimum value of $\log_{10}(f(\cdot) - \sigma)$ over the $N$ iterates.
        }
    \label{fig:subgradient}
\end{figure}

The performance of Algorithm \ref{alg-circumcenter} is shown in Figure \ref{fig:subgradient}.  The plot clearly shows the $O(\frac{1}{\sqrt{n}})$ behavior relative to the number of iterations $n$.  For comparison, in this example the set $A$ has diameter $\sqrt{10 + 4\sqrt{5}} < 5$, so an upper bound on the error in approximating the circumradius, by Corollary~\ref{circumcenter-complexity}, is $\frac{10}{\sqrt{n}}$.

\subsubsection*{Acknowledgements}
The implementation of the horospherical subgradient algorithm and the results in Figure \ref{fig:subgradient} are due to Ariel Goodwin.


\def\cprime{$'$} \def\cprime{$'$}

\end{document}